\documentclass[reqno,11pt]{amsart}

\usepackage{amsmath,amssymb,amsthm}
\usepackage{cases} 
\usepackage{tikz-cd}        %diagram

\usepackage[
backend=biber,
style=numeric,%alphabetic,
sorting=ynt
]{biblatex}
\usepackage{hyperref}
\hypersetup{colorlinks}
\hypersetup{
    bookmarksopen=true,     % show bookmarks bar?
    unicode=false,          % non-Latin characters in Acrobat’s bookmarks
    pdftoolbar=true,        % show Acrobat’s toolbar?
    pdfmenubar=true,        % show Acrobat’s menu?
    pdffitwindow=false,     % window fit to page when opened
    pdfstartview={FitH},    % fits the width of the page to the window
    pdftitle={My title},    % title
    pdfauthor={Author},     % author
    pdfsubject={Subject},   % subject of the document
    pdfcreator={Creator},   % creator of the document
    pdfproducer={Producer}, % producer of the document
    pdfkeywords={keyword1} {key2} {key3}, % list of keywords
    pdfnewwindow=true,      % links in new PDF window
    colorlinks=true,       % false: boxed links; true: colored links
    linkcolor=red,          % color of internal links (change box color with linkbordercolor)
    citecolor=magenta,        % color of links to bibliography
    filecolor=violet,      % color of file links
    urlcolor=blue      % color of external links
}
\newtheorem{theorem}{Theorem}[section]
\newtheorem{lemma}[theorem]{Lemma}
\newtheorem{corollary}{Corollary}[theorem]
\newtheorem{definition}[theorem]{Definition}
\newtheorem{proposition}[theorem]{Proposition}

\begin{document}
\sffamily %applying sans serif font

% 전문
\title{Leafwise cohomological expression of dynamical zeta functions on foliated dynamical systems}
\author{Junhyeong Kim}
\date{}

\address{
Graduate School of mathematics \endgraf 
Kyushu University \endgraf
744 Motooka, Nishi-ku, \endgraf
Fukuoka 819-0395, Japan
}

%\keywords{}

\thanks{$^*$Research Fellow of the Japan Society for the Promotion of Science.}

\email{j-kim@math.kyushu-u.ac.jp}


\begin{abstract}
A Riemmanian foliated dynamical system of 3-dimension $(\mathrm{RFDS}^{3})$ is a closed Riemannian 3-manifold with additional structures: foliation, dynamical system. 
In the context of arithmetic topology, it is a geometric/analytic analogue of an arithmetic scheme with a conjectural dynamical system suggested by C. Deninger.
In this paper, we show leafwise cohomological expression of dynamical zeta function on a Riemannian foliated dynamical system. 
\end{abstract}
\maketitle

% \setcounter{tocdepth}{1}
% \tableofcontents

% 서문
\section{Introduction}
% 연구 주제,필요성, 중요성

% 선행 연구, 밝히지 못한 점
In a series of papers (c.f. \cite{deninger1998some},\cite{deninger2002number},\cite{deninger2005arithmetic},\cite{deninger2006dynamical},\cite{deninger2008analogies}), C. Deninger considered arithmetic schemes $\overline{\mathrm{spec}\,\mathcal{O}_{K}}$ with a conjectural dynamical system for a number field $K/\mathbb{Q}$.
He interpreted the completed Dedekind zeta function $\hat{\zeta}_{K}(s)$ of $K$ in terms of infinite dimensional cohomology groups $H^{\bullet}_{\mathrm{dyn}}(\overline{\mathrm{spec}\,\mathcal{O}_{K}},\mathcal{R})$:
\begin{equation*}
    \hat{\zeta}_{K}(s)=\prod_{i=0}^{2} \mathrm{det}_{\infty} \left(\frac{1}{2\pi}(s-\Theta)|H^{i}_{\mathrm{dyn}}(\overline{\mathrm{spec}\,\mathcal{O}_{K}},\mathcal{R}) \right)^{(-1)^{i+1}},
\end{equation*}
where $\det_{\infty}$ denotes the zeta-regularized determinant and $\Theta$ denotes an infinitesimal generator of the flow.

This idea is extended to smooth closed 3-manifold $M$ with 1-codimensional foliation structure $\mathcal{F}$, transverse flow $\phi$ and a bundle-like metric $g_{\mathcal{F}}$ via Arithmetic topology (\cite{morishita2011knots}).
We call the manifold with the additional structure a \textbf{Riemannian foliated dynamical system} of 3 dimension, simply $\mathrm{RFDS}^3$.
It is a geometric/dynamical analogue of the above arithmetic scheme with a conjectural dynamical system, where closed orbits correspond to finite primes.
Note that dynamical zeta function corresponds to Dedekind zeta function in this context.
%논문의 목적, 방법,결과
The purpose of this paper is to show leafwise cohomological expression of dynamical zeta function on a Riemannian foliated dynamical system.
We describe our main results in the following. \\

Let $(M, \mathcal{F}, \phi, g_{\mathcal{F}})$ be a $\mathrm{RFDS}^3$.
The additional structures give a \textbf{reduced leafwise cohomology} $\bar{H}_{\mathcal{F}}^{\bullet}(M)$ and an \textbf{infinitesimal generator} $\Theta$.
Then we consider infinite series
\begin{equation*}
    \xi_{p}(s,z):=\sum_{\rho\in\mathrm{Sp}(\Theta_{p})}(s-\rho)^{-z},
\end{equation*}
where $\Theta_{p}$ denotes the operator $\Theta$ acting on $\bar{H}^{p}_{\mathcal{F}}(M)$.
We have our first main theorem:
\begin{theorem}\label{theorem:1.1}
The following assertions hold:
\begin{enumerate}
    \item {
        The series $\xi_{p}(s,z)$ is absolutely convergent on $\mathrm{Re}(z) >> 0$ for any $s\in\mathbb{C}$.
    }
    \item {
    It extends to a meromorphic function of $z\in\mathbb{C}$ and $s\in\mathbb{C}$} which is holomorphic at $z=0$.
\end{enumerate}
\end{theorem}
As a consequence of the theorem, the series defines a Hurwitz-type spectral zeta function associated with the infinitesimal generator.
It follows that a zeta-regularized determinant can be defined since the spectral zeta function is regular at $z=0$.

For the second result, we define the dynamical zeta function for $\mathrm{RFDS}^3$.
\begin{equation*}
    \zeta_{\mathcal{F}}(s) = \prod_{\gamma}(1-e^{-s\cdot{l(\gamma)}})^{-\epsilon_{\gamma}},
\end{equation*}
where $\gamma$ runs over closed orbits of $\phi$ and $l(\gamma)$ is the length of $\gamma$.
Here, $\epsilon_{\gamma}$ denotes the index of a closed orbit.
We give our second main theorem as follows:
\begin{theorem}\label{theorem:1.2}
The dynamical zeta function on a Riemmanian foliated dynamical system of 3 dimension has a leafwise cohological expression
\begin{equation*}
    \zeta_{\mathcal{F}}(s) = \prod_{i=0}^{2} \mathrm{det}_{\infty}(s-\Theta | \bar{H}_{\mathcal{F}}^{i}(M))^{(-1)^{i+1}}.
\end{equation*}
\end{theorem}

%논문의구성 순서
The contents of this paper are organized as follows:
In section 2, 3, 4, we introduce a Riemannian foliated dynamical system of 3 dimension ($\mathrm{RFDS}^3$) and basic notions: leafwise cohomology and infinitesimal generator.
In section 5, we give a proof of the main theorem \ref{theorem:1.1}.
In section 6, 7, we recall the zeta-regularized determinant and dynamical zeta function for $\mathrm{RFDS}^3$.
In section 8, we show a leafwise cohomlogical expression of the dynamical zeta function on $\mathrm{RFDS}^3$.

% part 1
\section{Riemannian foliated dynamical system (RFDS\texorpdfstring{$^{3}$}{3})}

We consider a smooth, compact, orientable, closed 3-manifold $M$ with additional structure: foliation $\mathcal{F}$, transverse flow $\phi$. 

    \subsection{Foliation}
    A foliation $\mathcal{F}$ of $d$-codimension is a partition of sub-manifolds of $d$-codimension.
    Let $M$ be a smooth, connected, closed and oriented manifold of $n$-dimension. It is equipped with a foliation of $d$-codimension as follows: let $(U_{i}, \phi_{i})_{i\in{I}}$ be an atlas.
    The transition maps $\phi_{ij}:=\phi_{j}\circ\phi_{i}^{-1}$ which are defined over $U_{i}\cap{U_{j}}$ take forms 
    \begin{multline*}
        \phi_{ij}(x_{1},\cdots,x_{d},y_{d+1},\cdots,y_{n}) = (\phi_{ij}^{1}(x_{1},\cdots,x_{d}),\cdots,\phi_{ij}^{d}(x_{1},\cdots,x_{d}), \\
        \phi_{ij}^{d+1}(x_{1},\cdots,x_{d},y_{d+1},\cdots,y_{n}), \phi_{ij}^{n}(x_{1},\cdots,x_{d},y_{d+1},\cdots,y_{n}))
    \end{multline*}
    for $i,j \in {I}$.
    By piecing together the stripes, where $(x_{1}, \cdots, x_{d})$ are constant, from chart to chart, we obtain a maximal immersed sub-manifold $\mathcal{L}$ whose first $d$ local coordinates are constant on each $U_{i}$.
    We call the sub-manifold a $\textit{leaf}$ of the foliation.
    The foliation consists of the disjoint union of leaves.
    
    \subsection{Transversal flow}
    A transverse flow $\phi$ is a smooth $\mathbb{R}$-action on a manifold $M$
    \begin{equation*}
        \phi : \mathbb{R} \times {M} \rightarrow {M}
    \end{equation*}
    which maps leaves of a foliation to leaves. 
    For any two points $x$ and $y$ in a same leaf $\mathcal{L}$, there is a leaf $\mathcal{L}^{'}$ containing $\phi(t,x)$ and $\phi(t, y)$ for any $t\in\mathbb{R}$.
    Let $\dot{\phi}$ be the vector field giving the velocity vector at a point.
    We denote by $\omega_{\phi}$ the dual 1-form of the vector field $\dot{\phi}$ .
    
    \subsection{Bundle-like metric}
    A Riemannian metric $g_{\mathcal{F}}$ on $(M,\mathcal{F},\phi)$ is called a \textbf{bundle-like metric} whose geodesics are perpendicular to all leaves whenever they are perpendicular to one leaf.
    Note that any 1-codimensional foliation without singularities is Riemannian.\\

We consider 3-manifolds with additional structures as follows:
\begin{definition}
We define a foliated dynamical system on a 3-manifold by a triple $(M,\mathcal{F},\phi)$, where
\begin{enumerate}
    \item $M$ is a smooth, compact, orientable 3-manifold,
    \item $\mathcal{F}$ is a 1-codimensional foliation on $M$,
    \item $\phi$ is a smooth $\mathbb{R}$-action acting on $M$ such that
    \begin{enumerate}
        \item The flow is transverse to the leaves of the foliation up to a finite number of compact leaves;
        \item The $\mathbb{R}$-action maps leaves to leaves.
    \end{enumerate}
\end{enumerate}
\end{definition}
The manifold and the flow may have boundaries and fixed-points.
For this paper, we assume that the manifold is closed and the flow has no fixed-point.
It is known that only mapping torus allows such a foliated dynamical system on itself.
\section{Leafwise cohomology}

    \subsection{Leafwise de Rham complex}
    
    For the triple $(M,\mathcal{F},\phi)$, let $T\mathcal{F}$ be a sub-bundle of the tangent bundle $TM$ which is tangent to the leaves of the foliation.
    The restriction of $T\mathcal{F}$ on a leaf $\mathcal{L}$ is identified with the tangent bundle $T\mathcal{L}$ of the leaf.
    We call $T\mathcal{F}$ the \textit{leafwise tangent bundle}.

    We define the space of leafwise $i$-forms by
    \begin{equation*}
        \Omega^{i}_{\mathcal{F}}(M):=\Gamma(M,\wedge^{i}T^{*}\mathcal{F})\subset{\Omega^{i}(M)}.
    \end{equation*}
    Let $d_{\mathcal{F}}$ (resp. $d_{0}$) be the exterior derivative acting only along leaves (resp. the flow).
    Then the de Rham complex $(\Omega^{i}(M),d^{i})$  has a decomposition:
    \begin{center}
    \begin{tikzcd}
    \cdots \ar[r]\ar[dr, ""] & \Omega^{i}_{\mathcal{F}}(M) \ar[d, "\oplus"]\ar[r, "d_{\mathcal{F}}^{i}"]\ar[dr, "d_{0}^{i}"] & \Omega^{i+1}_{\mathcal{F}}(M)  \ar[d, "\oplus"]\ar[r]\ar[dr, ""] & \cdots\\
    \cdots \ar[r] & \Omega^{i}_{0}(M) \ar[r, "d_{\mathcal{F}}^{i}"] & \Omega^{i+1}_{0}(M) \ar[r] & \cdots
    \end{tikzcd}
    \end{center}
    where $\Omega_{0}^{i}(M)$ is the complement of $\Omega^{i}_{\mathcal{F}}(M)$.
    
    We simply denote the restriction $d^{i}|_{\Omega^{i}_{\mathcal{F}}(M)}$ by $d_\mathcal{F}^{i}$.
    Since we have $d^{i+1}_{\mathcal{F}}\circ{d}^{i}_{\mathcal{F}}=0$ on $\Omega^{i}_{\mathcal{F}}(M)$, 
    the pairs $\{(\Omega^{i}_{\mathcal{F}}(M), d^{i}_{\mathcal{F}})\}_{i}$ form a cochain complex:
    \begin{equation*}
        0\overset{}{\rightarrow}\Omega^{0}_{\mathcal{F}}(M)\overset{d^{0}_{\mathcal{F}}}{\rightarrow}\Omega^{1}_{\mathcal{F}}(M)\overset{d^{1}_{\mathcal{F}}}{\rightarrow}\Omega^{2}_{\mathcal{F}}(M)\overset{d^{2}_{\mathcal{F}}}{\rightarrow}0.
    \end{equation*}
    We call the complex \textit{leafwise de Rham complex}.
    
    \subsection{Leafwise cohomology}
    We denote the kernel of $d_{\mathcal{F}}^{i}$ by $Z_{\mathcal{F}}^{i}(M)$ and the image of $d_{\mathcal{F}}^{i}$ by $B_{\mathcal{F}}^{i+1}(M)$.
    Note that a leafwise $i$-th form in $Z_{\mathcal{F}}^{i}(M)$ (resp. $B_{\mathcal{F}}^{i+1}(M)$) is called leafwise closed $i$-th form (resp. leafwise exact $i$-th form).
    
    \begin{definition}[Leafwise cohomology]
    We define the $i$-th leafwise cohomology group by
    \begin{equation*}
        H^{i}_{\mathcal{F}}(M):={Z_{\mathcal{F}}^{i}(M)}/{B_{\mathcal{F}}^{i}(M)}.
    \end{equation*}
    \end{definition}
    The leafwise cohomology group is trivial for $i>2$.
    
    Unfortunately, the leafwise cohomology group is of infinite dimension in general and not a Hausdorff space.
    We modify it by taking a quotient with respect to the closure in the smooth topology.
    \begin{equation*}
        \bar{H}^{i}_{\mathcal{F}}(M):= {Z_{\mathcal{F}}^{i}(M)}/{\overline{B_{\mathcal{F}}^{i}(M)}}.
    \end{equation*}
    We call it the reduced leafwise cohomology group.
    
    % \subsection{Scalar product}
    % Let $<,>_{\mathcal{F}}$ be the induced metric on $\wedge^{i}T^{*}\mathcal{F}$ from the bundle-like metric $g_{\mathcal{F}}$.
    % If we choose an orientation on $M$, it determines a volume form $\mathrm{vol}$ of $M$.
    % We set a scalar product on $\Omega^{i}_{\mathcal{F}}(M)$ by
    % \begin{equation*}
    %     (\omega,\eta)_{\mathcal{F}} = \int_{M}<\omega,\eta>_{\mathcal{F},x}\mathrm{vol}(x),
    % \end{equation*}
    % where $\omega$, $\eta$ are leafwise $i$-th forms.
    % The scalar product induces $\bar{H}^{i}_{\mathcal{F}}(M)$ a Hilbert structure.
    
    \subsection{Leafwise Hodge theorem}
    
    For a bundle-like metric $g_{\mathcal{F}}$, we have the Hodge $*$-operator.
    If we denote by $\delta$ the adjoint operator of the exterior derivative $d$, it has a decomposition into $\delta_{\mathcal{F}}$ and $\delta_{0}$
    \begin{center}
    \begin{tikzcd}
    \cdots  & \ar[l] \Omega^{i}_{\mathcal{F}}(M) \ar[d, "\oplus"] & \ar[l, "\delta_{\mathcal{F}}^{i}"] \Omega^{i+1}_{\mathcal{F}}(M)  \ar[d, "\oplus"] & \ar[l] \cdots\\
    \cdots & \ar[l]\ar[ul, ""] \Omega^{i}_{0}(M)  & \ar[l, "\delta_{\mathcal{F}}^{i}"]\ar[ul, "\delta_{0}^{i}"] \Omega^{i+1}_{0}(M) & \ar[l]\ar[ul, ""] \cdots.
    \end{tikzcd}
    \end{center}
    
    \begin{definition}[Leafwise Laplacian]
    An operator defined by
    \begin{equation*}
        \Delta_{\mathcal{F}} := d_{\mathcal{F}}\delta_{\mathcal{F}} + \delta_{\mathcal{F}}d_{\mathcal{F}} \text{ on  } \Omega_{\mathcal{F}}^{i}(M)
    \end{equation*}
    is called the \textbf{leafwise Laplacian}.
    A leafwise form $\omega\in\ker\Delta_{\mathcal{F}}$ is called a leafwise harmonic form.
    \end{definition}
    
    If we define $\Delta_{0}$ on $\Omega_{\mathcal{F}}^{i}(M)$ by $\delta_{0}d_{0}$, the restriction of the Laplacian $\Delta$ on $\Omega_{\mathcal{F}}^{i}(M)$ can be represented by
    \begin{equation*}
        \Delta|_{\Omega_{\mathcal{F}}^{i}(M)} = \Delta_{\mathcal{F}} + \Delta_{0}.
    \end{equation*}

    We have a significant proposition by \'Alvarez L\'opez and Kordyukov (\cite{lopez2001long}):
    \begin{proposition}[Leafwise Hodge theorem]
    Given a bundle-like metric, An leafwise cohomology class can be uniquely represented by a leafwise harmonic form.
    We have an isomorphism
    \begin{equation*}
        \bar{H}_{\mathcal{F}}^{i}(M) \cong \ker(\Delta_{\mathcal{F}}^{i}).
    \end{equation*}
    \end{proposition}
\section{Infinitesimal generator}

Assume that the $\mathbb{R}$-action $\phi$ is conformal on $\Omega_{\mathcal{F}}^{i}(M)$ with respect to the bundle-like metric $g_{\mathcal{F}}$, i.e.
\begin{equation*}
    (\phi^{t*}\omega,\phi^{t*}\eta)=(\omega,\eta)\text{ for } \forall{t}\in\mathbb{R}.
\end{equation*}
It is easy to check that $\phi^{t*}$ on $\bar{H}^{i}_{\mathcal{F}}(M)$ is surjective and strongly continuous, i.e.
\begin{equation*}
    \lim_{t\rightarrow{t_{0}}}\phi^{t*}[h]=\phi^{t_{0}*}[h]\quad\text{for}\,\,\forall{t}_{0}\in\mathbb{R}, [h]\in\bar{H}^{i}_{\mathcal{F}}(M).
\end{equation*}
Then the following lemma follows from the Stone's theorem
\begin{lemma}[Stone's theorem]\label{lem:1.3}
We define the infinitesimal generator of $(\phi^{t*})_{t\in\mathbb{R}}$ by
\begin{equation*}
    \Theta := \underset{t\rightarrow0}{\mathrm{lim}} \frac{\phi^{t*}- \mathrm{id}}{t}.
\end{equation*}
Since $(\phi^{t*})_{t\in\mathbb{R}}$ is the strongly continuous one-parameter unitary group on the Hilbert space $\bar{H}^{i}_{\mathcal{F}}(M)$,
then $A:=-i\Theta$ is self-adjoint on $\bar{H}^{i}_{\mathcal{F}}(M)$ and we have
\begin{equation*}
    \phi^{t*}=e^{itA}=e^{t\Theta}\text{ for }\forall{t}\in\mathbb{R}.
\end{equation*}
\end{lemma}
The infinitesimal generator is a first-order differential operator along the transversal flow.
Then the exterior derivative $d_{0}$ along the flow (resp. adjoint operator $\delta_{0}$) can be represented by
\begin{equation*}
\begin{split}
    & d_{0}\omega = \Theta\omega \wedge \omega_{\phi},\\
    & \delta_{0}(\omega \wedge \omega_{\phi}) = -\Theta\omega.
\end{split}
\end{equation*}
It leads the following lemma.
\begin{lemma}\label{lem:5.2}
The negative square of the infinitesimal generator on a leafwise cohomology group coincide with the Laplacian on a space of leafwise harmonic forms
\begin{equation*}
    \begin{split}
        -\Theta^{2}|_{\bar{H}_{\mathcal{F}}^{i}(M)} &= \Delta_{0}|_{\ker(\Delta_{\mathcal{F}})} \\
        &= \Delta|_{\ker(\Delta_{\mathcal{F}})}.
    \end{split}
\end{equation*}
\end{lemma}
Since $M$ is compact and closed, the Laplacian has pure point spectrum which consists of non-negative eigenvalues with finite multiplicity.
Hence the infinitesimal generator has pure imaginary eigenvalues with finite multiplicity.\\

For a leafwise harmonic form $\omega_{\mathcal{F}}\in\ker\Delta_{\mathcal{F}}$, a fundamental solution of the heat equation whose initial value is $\omega_{\mathcal{F}}$, i.e.
\begin{equation*}
    \begin{cases}
        \left(\frac{\partial}{\partial{t}} + \Delta_{0} \right) \omega_{\mathcal{F}}^{t} = 0, \\
        \omega_{\mathcal{F}}^{0} = \omega_{\mathcal{F}},
    \end{cases}
\end{equation*}
is given by
\begin{equation*}
        \omega_{\mathcal{F}}^{t}(x,y,s) =
        \begin{cases}
            \int_{S^{1}} K(t,s,s')\omega_{\mathcal{F}}(x,y,s')ds' & \text{if } (x,y,s) \text{ is periodic},\\
            \int_{\mathbb{R}} K(t,s,s')\omega_{\mathcal{F}}(x,y,s')ds' & \text{otherwise}
        \end{cases}
\end{equation*}
where $K(t,s,s')$ is a factor of the heat kernel of the Laplacian $\Delta$.
Hence we have an asymptotic expansion for the heat kernel around $t=0$
\begin{equation*}
    \mathrm{tr}(e^{-\Delta_{0}t}|\ker\Delta_{\mathcal{F}}) \overset{t\downarrow0}{\sim} t^{-\frac{1}{2}}(a_{0} + a_{1}t + a_{2}t^{2} + \cdots).
\end{equation*}
Then the spectral zeta function $\zeta_{\Delta_{0}}(s)$ associated with $\Delta_{0}$ has only simple poles at $s=\frac{1}{2}-n$ $(n=0,1,2, \cdots)$.

Fixing a positive number $T>0$, we consider a series
\begin{equation*}
    V^{p}(t) = \sum_{\text{Im}(\rho)>T} e^{\rho{it}}
\end{equation*}
where $\rho$ runs over the spectrum of $\Theta$ on $\bar{H}_{\mathcal{F}}^{p}(M)$.
It follows from the lemma \ref{lem:5.2} that it is a partial sum of $\mathrm{tr}(e^{-\sqrt{\Delta_{0}}t}|\ker\Delta_{\mathcal{F}})$.
Since $\mathrm{tr}(e^{-\sqrt{\Delta_{0}}t}|\ker\Delta_{\mathcal{F}})$ is the inverse Mellin transform of $\Gamma(s)\zeta_{\Delta_{0}}(\frac{s}{2})$ with simple poles at $s=1, -2n$  and double poles at $s=-2n-1$ $(n=0,1,2,\cdots)$, we deduce that the series $V^{p}(t)$ converges absolutely and has an asymptotic expansion around $t=0$ as follows:
\begin{equation*}
    V^{p}(t) \overset{t\downarrow0}{\sim} at^{-1} + \sum_{k=0}^{N} (b_{k} + c_{k}{t}\log{t})t^{2k} + \mathcal{O}_{1}(t^{N-1}) + \mathcal{O}_{2}(t^{N-1})t\log{t}.
\end{equation*}
\section{Proof of theorem \ref{theorem:1.1}}

\begin{proof}
We fix a positive number $T>0$.
We consider the 2 series for $s\in\mathbb{C}$ such that $|\mathrm{Im}(s)|<T$
\begin{equation*}
    \begin{split}
        \theta_{p}^{+}(t) &= V^{p}(t) e^{-sit}, \\
        \theta_{p}^{-}(t) &= V^{p}(t) e^{sit}.
    \end{split}
\end{equation*}
They play a role like a partition function.

We take the Mellin transform for the series and define the following functions
\begin{equation*}
    \begin{split}
        \xi_{p}^{+}(s,z) &= \frac{e^{\frac{\pi}{2}iz}}{\Gamma(z)} \int_{0}^{\infty}\theta_{p}^{+}(t)t^{z-1}dt,\\
        \xi_{p}^{-}(s,z) &= \frac{e^{-\frac{\pi}{2}iz}}{\Gamma(z)}  \int_{0}^{\infty}\theta_{p}^{-}(t)t^{z-1}dt.
    \end{split}
\end{equation*}

Since $V^{p}(t)$ is convergent, we have for $\mathrm{Re}(z) > 1$
\begin{equation*}
    \begin{split}
        \xi_{p}^{+}(s,z) &= \sum_{\mathrm{Im}(\rho)>T}(s-\rho)^{-z}, \\
        \xi_{p}^{-}(s,z) &= \sum_{\mathrm{Im}(\rho)<-T}(s-\rho)^{-z}.
    \end{split}
\end{equation*}

Next, we consider
\begin{equation*}
    \begin{split}
        \xi_{p}^{+}(s,z) &= \frac{e^{\frac{\pi}{2}iz}}{\Gamma(z)} \int_{0}^{\infty}\theta_{p}^{+}(t)t^{z-1}dt \\
        &= \frac{e^{\frac{\pi}{2}iz}}{\Gamma(z)} \left( \int_{0}^{1}\theta_{p}^{+}(t)t^{z-1}dt + \int_{1}^{\infty}\theta_{p}^{+}(t)t^{z-1}dt \right), \\
        \xi_{p}^{-}(s,z) &= \frac{e^{-\frac{\pi}{2}iz}}{\Gamma(z)} \int_{0}^{\infty}\theta_{p}^{-}(t)t^{z-1}dt \\
        &= \frac{e^{-\frac{\pi}{2}iz}}{\Gamma(z)} \left( \int_{0}^{1}\theta_{p}^{-}(t)t^{z-1}dt + \int_{1}^{\infty}\theta_{p}^{-}(t)t^{z-1}dt \right).
    \end{split}
\end{equation*}
Since $V^{p}(t)$ is of rapid decay at infinity, the second terms are convergent for any $z\in\mathbb{C}$.
Hence we have
\begin{equation*}
    \begin{split}
        \xi_{p}^{+}(s,z) &= \frac{e^{\frac{\pi}{2}iz}}{\Gamma(z)} \left( \int_{0}^{1}\theta_{p}^{+}(t)t^{z-1}dt + \int_{1}^{\infty}\theta_{p}^{+}(t)t^{z-1}dt \right) \\
        &= \frac{e^{\frac{\pi}{2}iz}}{\Gamma(z)} \left ( \frac{a}{z-1} - \frac{asi + b_{0}}{z} - \frac{b_{0}si}{z+1} - \frac{c_{0}}{(z+1)^{2}} + \cdots \right) \\
        &= \frac{ai}{z-1} + \eta_{p}^{+}(s,z),
    \end{split}
\end{equation*}
where $\eta_{p}^{+}(s,z)$ is a meromorphic function of $(s,z)$ for $|\mathrm{Im}(s)|<T$ and $z\in\mathbb{C}$ and is regular at $z=0$.
The meromorphic function $\eta_{p}^{+}(s,z)$ has only simple poles at $z=-2n-1$ $(n=0,1,2,\cdots)$. 
% It follows that the series $\sum_{\mathrm{Im}(\rho)>T}(s-\rho)^{-z}$ has a meromorphic continuation to all of $\mathbb{C}$ and is regular at $z=0$.
The same result holds for $\xi_{p}^{-}(s,z)$ by the same argument:
\begin{equation*}
    \xi_{p}^{-}(s,z) = \frac{-ai}{z-1} + \eta_{p}^{-}(s,z).
\end{equation*}
Note that $\eta_{p}^{-}(s,z)$ is meromorphic in $|\mathrm{Im}(s)|<T$ and $z\in\mathbb{C}$ and is regular at $z=0$.
Since $\xi_{p}(s,z)$ differs from $\xi_{p}^{+}(s,z)+\xi_{p}^{-}(s,z)$ by the sum of the finite terms, we have that $\xi_{p}(s,z)$ is a meromorphic function for all $|\mathrm{Im}(s)|<T$ and all $z\in\mathbb{C}$ and is regular at $z=0$.
\end{proof}

\section{Zeta-regularized determinant}

We recall the notion of the zeta-regularized determinant.
Let $\Theta:V\rightarrow{V}$ be a linear operator acting on a complex vector space $V$ of countable dimension.
We assume that $V$ is the direct sum of finite-dimensional $\Theta$-invariant sub-spaces.
Let $\mathrm{Sp}(\Theta)$ be the set of eigenvalues of $\Theta$. 
The spectral zeta function associated with the operator $\Theta$ is defined by the analytic continuation of Dirichlet series
\begin{equation*}
  \zeta_{\Theta}(s)=\sum_{\lambda\ne0\in\mathrm{Sp}(\Theta)}\lambda^{-s}\text{ with }\lambda^{-s}=|\lambda|^{-s}e^{-is(\mathrm{Arg}\lambda)}, -\pi<\mathrm{Arg}\lambda\le\pi.
\end{equation*}
We assume that the Dirichlet series converges absolutely on some right-half plane and has an analytic continuation to the half plane $\mathrm{Re}(s)>-\epsilon$ for some $\epsilon>0$ which is holomorphic at $s=0$.
Under these conditions, we define a zeta-regularized determinant by
\begin{equation*}
    \mathrm{det}_{\infty}(\Theta|V):=\exp\left(-\partial_{s}\zeta_{\Theta}(0)\right).
\end{equation*}

\section{Dynamical zeta function on RFDS\texorpdfstring{$^{3}$}{3}}

Let $(M,\mathcal{F},\phi,g_{\mathcal{F}})$ be the foliated dynamical system with a bundle-like metric which we discussed above.
We define the dynamical zeta function for a $\mathrm{RFDS}^{3}$ by the analytic continuation of the infinite product
\begin{equation*}
  \zeta_{\mathcal{F}}(s)=\prod_{\gamma}(1-e^{-s\cdot l(\gamma)})^{-\epsilon_{\gamma}},
\end{equation*}
where $\gamma$ runs over periodic orbits of $\phi$ and  $l(\gamma)$ is the length of $\gamma$. 
Here, $\epsilon_{\gamma}$ is the index of a closed orbit.

    \subsection{Index of a closed orbit}
    For a closed orbit $\gamma$ of $\phi$, we set an index
    \begin{equation*}
        \epsilon_{\gamma}:=\mathrm{sgn}\,\det(1-T_{x}\phi^{l(\gamma)}|T_{x}\mathcal{F})\quad\mathrm{for}\,{x}\in\gamma,
    \end{equation*}
    where $T_{x}\phi^{t}:T_{x}\mathcal{F}\rightarrow{T_{\phi^{t}(x)}\mathcal{F}}$ is the differential of $\phi^{t}$.
    It does not depend on the choice of the point $x\in\gamma$.
    We call a closed orbit $\gamma$ \textbf{non-degenerate} in a sense that $\epsilon_{\gamma}$ is non-zero.

    \subsection{Absolute convergent condition}
    It is known that the infinite product converges absolutely on $\mathrm{Re}(s)>h(\phi)$ where $h(\phi)$ is \textbf{topological entropy} and only if it is finite. 
    Note that  the topological entropy $h(\phi)$ is defined by
    \begin{equation*}
        h(\phi):=\lim_{T\rightarrow+\infty}\frac{1}{T}\log{N(T)}\ge0,
    \end{equation*}
    where  $N(T)$ denotes the cardinality of orbits whose length is less than or equal to $T$, i.e. $N(T)=\mathrm{Card}\{\gamma|l(\gamma)\le{T}\}$. 
    We assume that the topological entropy $h(\phi)$ of a foliated dynamical system  $(M,\mathcal{F},\phi)$ is finite so that $\zeta(s)$ converges absolutely on the right-half plane.
    
\section{Proof of theorem \ref{theorem:1.2}}

    \subsection{Dynamical Lefschetz trace formula}
    For the foliated dynamical system $(M,\mathcal{F},\phi, g_{\mathcal{F}})$ whose closed orbits are all non-degenerate, Alvarez-Lopez and Kordyukov developed the \textbf{dynamical Lefschetz trace formula}:
    \begin{proposition}[\cite{lopez2002distributional}]
    For every test function $\varphi\in\mathcal{D}(\mathbb{R})=C^{\infty}_{0}(\mathbb{R})$, the operator
    \begin{equation*}
        A_{\varphi}=\int_{\mathbb{R}}\varphi(t)\phi^{t*}dt
    \end{equation*}
    on $\bar{H}^{i}_{\mathcal{F}}(M)$ is of trace class.
    Setting:
    \begin{equation*}
        \mathrm{Tr}(\phi^{t*}|\bar{H}^{i}_{\mathcal{F}}(M))
        =\mathrm{tr}A_{\varphi}
    \end{equation*}
    defines a distribution on $\mathbb{R}$.
    The following formula holds in $\mathcal{D}^{'}(\mathbb{R})$:
    \begin{equation*}
        \sum_{i=0}^{\dim{\mathcal{F}}}(-1)^{i}\mathrm{Tr}(\phi^{t*}|\bar{H}^{i}_{\mathcal{F}}(M))
        =\chi_{\mathrm{Co}}(\mathcal{F},g_{\mathcal{F}})\delta_{0}
        +\sum_{\gamma}l(\gamma)\sum_{k\in\mathbb{Z}\backslash0}\epsilon_{\gamma}\delta_{kl(\gamma)}.
    \end{equation*}
    Here $\chi_{\mathrm{Co}}(\mathcal{F},\mu)$ denotes Connes' Euler characteristic of the foliation with respect to the bundle-like metric (c.f. \cite{connes2000noncommutative}) and $\delta_{\tau}$ is the Dirac delta function in $\mathcal{D}^{'}(\mathbb{R})$ which is non-zero at $\tau$.
    \end{proposition}
    
    The lemma \ref{lem:1.3} (Stone's theorem) leads to the corollary:
    \begin{corollary} \label{cor:2.1.1.}
    The following equality holds in $\mathcal{D}^{'}(\mathbb{R}_{>0})$:
    \begin{equation*}
        \sum_{i=0}^{2}(-1)^{i}\mathrm{Tr}(\phi^{t*}|\bar{H}^{i}_{\mathcal{F}}(M))
        =\sum_{i=0}^{2}(-1)^{i}\sum_{\rho\in\mathrm{Sp}(\Theta_{i})}e^{\rho{t}}.
    \end{equation*}
    where $\Theta_{i}$ denotes the operator $\Theta$ acting on $\bar{H}^{i}_{\mathcal{F}}(M)$.
    \end{corollary}

It is enough to show that the zeta-regularized determinant coincides with the infinite product of periodic orbits on some right-half plane of the topological entropy, $i.e.$,
\begin{equation*}
    \begin{split}
        \zeta_{\mathcal{F}}(s) 
        &= \prod_{\gamma}(1-e^{-sl(\gamma)})^{-\epsilon_{\gamma}}\\
        &= \prod_{i=0}^{2}\det_{\infty}(s-\Theta|\bar{H}^{i}_{\mathcal{F}}(M))^{(-1)^{i+1}} \quad\text{for }\mathrm{Re}(s)>{P},
    \end{split}
\end{equation*}
where $P$ is a sufficiently large number $P>h(\phi)$.
Then the assertion follows from the uniqueness of analytic continuation.

We apply the Laplace transform for the dynamical Lefschetz trace formula.
We have
\begin{equation}\label{lhs}
    \mathcal{L}\left[ \sum_{i=0}^{2}(-1)^{i}\sum_{\rho\in\mathrm{Sp}(\Theta_{i})} t^{z-1}e^{\rho{t}} \right](s) = \Gamma(z)\sum_{i=0}^{2}(-1)^{i}\sum_{\rho\in\mathrm{Sp}(\Theta_{i})}(s-\rho)^{-z}
\end{equation}
for the left hand side, and
\begin{equation}\label{rhs}
    \mathcal{L}\left[ \sum_{\gamma}\sum_{n\in\mathbb{N}}l(\gamma)\epsilon_{\gamma}\delta_{nl(\gamma)}t^{z-1} \right](s) 
    = \sum_{\gamma}\sum_{n\in\mathbb{N}}\frac{l(\gamma)\epsilon_{\gamma}}{e^{nl(\gamma)s}}(nl(\gamma))^{z-1}
\end{equation}
for the right hand side.
Both sides are defined for $\mathrm{Re}(z)>1$ over where the former infinite series (\ref{lhs}) is defined from the proof of theorem \ref{theorem:1.1}, and $\mathrm{Re}(s)> h(\phi)$ over where the latter infinite series (\ref{rhs}) is defined.
We denote by $P$ a sufficiently large number bigger than $h(\phi)$.

Let $\mathrm{L}_{\delta-}$ be a contour consisting of the lower edge of the cut from $-\infty$ to $-\delta$, the circle $t=\delta{e^{i\phi}}$ for $-\pi\le\phi\le\pi$ and the upper edge of the cut from $-\delta$ to $-\infty$.
\begin{equation*}
    \int_{\mathrm{L}_{\delta-}}e^{\lambda{t}}t^{-z}dt = 2i\sin(z\pi)\int_{\delta}^{\infty}e^{-v}v^{-z}dv + I
\end{equation*}
where $I$ denotes the integral along the circle $t=|\delta|$.
Since $I$ tends to zero as $\delta\rightarrow{0}$, we have
\begin{equation*}
    \begin{split}
        \lim_{\delta\rightarrow{0}}\int_{\mathrm{L}_{\delta-}}e^{\lambda{t}}t^{-z}dt 
        &= 2i\sin(z\pi)\Gamma(1-z)\\
        &= \frac{2\pi{i}}{\Gamma(z)}
    \end{split}
\end{equation*}
Hence we have the formula for $\lambda>0$
\begin{equation*}
    \frac{\lambda^{z-1}}{\Gamma(z)}=\frac{1}{2\pi{i}}\lim_{\delta\rightarrow{0}}\int_{\mathrm{L}_{\delta-}}e^{\lambda{t}}t^{-z}dt.
\end{equation*}
By applying the formula for the series (\ref{rhs}), we get
\begin{equation*}
    \begin{split}
        \frac{1}{\Gamma(z)}
        \sum_{\gamma}\sum_{n\in\mathbb{N}}
        \frac{l(\gamma)\epsilon_{\gamma}}{e^{nl(\gamma)s}}
        (nl(\gamma))^{z-1}
        &=\frac{1}{2\pi{i}}
        \lim_{\delta\rightarrow{0}}\int_{\mathrm{L}_{\delta-}}\left(\sum_{\gamma}\sum_{n\in\mathbb{N}} l(\gamma)\epsilon_{\gamma}e^{-nl(\gamma)(s-t)} \right)t^{-z}dt\\
        &=\frac{-1}{2\pi{i}}
        \lim_{\delta\rightarrow{0}}\int_{\mathrm{L}_{\delta-}}
        \frac{\zeta^{'}_{\mathcal{F}}}{\zeta_{\mathcal{F}}}(s-t)t^{-z}dt.
    \end{split}
\end{equation*}

Since the series (\ref{lhs}) has a meromorphic extension and is holomorphic at $z=0$ from theorem \ref{theorem:1.1}, we obtain the two equalities for $-\pi\le\mathrm{arg}(t)\le\pi$
\begin{equation*}
    \begin{split}
        \sum_{i=0}^{2}(-1)^{i}\xi_{i}(s,z) 
        &=\frac{-1}{2\pi{i}}
        \lim_{\delta\rightarrow{0}}\int_{\mathrm{L}_{\delta-}}\frac{\zeta^{'}_{\mathcal{F}}}{\zeta_{\mathcal{F}}}(s-t)t^{-z}dt, \\
        \sum_{i=0}^{2}(-1)^{i}\partial_{z}\xi_{i}(s,0) 
        &= \frac{1}{2\pi{i}}
        \lim_{\delta\rightarrow{0}}\int_{\mathrm{L}_{\delta-}}\frac{\zeta^{'}_{\mathcal{F}}}{\zeta_{\mathcal{F}}}(s-t)\log(|t|e^{\mathrm{arg}(t)i})dt
    \end{split}
\end{equation*}

It remains to see the following:
\begin{equation*}
    \begin{split}
        &\frac{1}{2\pi{i}}
        \lim_{\delta\rightarrow{0}}\int_{\mathrm{L}_{\delta-}}\frac{\zeta^{'}_{\mathcal{F}}}{\zeta_{\mathcal{F}}}(s-t)\log(|t|e^{\mathrm{arg}(t)i})dt\\
        &=\frac{1}{2\pi{i}}
        \int_{-\infty}^{0}\frac{\zeta^{'}_{\mathcal{F}}}{\zeta_{\mathcal{F}}}(s-t)
        ( \log(|t|) - \pi{i} )dt 
        +\frac{1}{2\pi{i}}
        \int_{0}^{-\infty}\frac{\zeta^{'}_{\mathcal{F}}}{\zeta_{\mathcal{F}}}(s-t)
        ( \log(|t|) + \pi{i} )dt\\
        &=\int_{0}^{-\infty}\frac{\zeta^{'}_{\mathcal{F}}}{\zeta_{\mathcal{F}}}(s-t)dt = -\int_{0}^{\infty}\frac{\zeta^{'}_{\mathcal{F}}}{\zeta_{\mathcal{F}}}(s+t)dt \\
        &=\log(\zeta_{\mathcal{F}}(s))
    \end{split}
\end{equation*}

Therefore we have
\begin{equation*}
    \begin{split}
        \zeta_{\mathcal{F}}(s) &= \prod_{i=0}^{2}\exp(-\partial_{z}\xi(s,0))^{(-1)^{i+1}} \\
        &= \prod_{i=0}^{2} \mathrm{det}_{\infty}(s-\Theta | \bar{H}^{i}_{\mathcal{F}}(M))^{(-1)^{i+1}}
    \end{split}
\end{equation*}
for $\mathrm{Re}(s)>P\ge{h(\phi)}$.
Hence the theorem \ref{theorem:1.2} follows from the uniqueness of  the analytic continuation.

% 참고문헌
\medskip
\printbibliography

% 끝
\end{document}